\documentclass[reqno,10pt]{amsart}
\usepackage{amsfonts}

\usepackage{graphicx}
\usepackage{pstricks}
\usepackage{amsmath}
\usepackage{amsmath}
\usepackage{amsxtra}

\newtheorem{theorem}{Theorem}

\newtheorem{corollary}[theorem]{Corollary}

\newtheorem{lemma}[theorem]{Lemma}

\newtheorem{proposition}[theorem]{Proposition}
\newtheorem{remark}[theorem]{Remark}

\begin{document}
\title{Quartically hyponormal weighted shifts \\
need not be $3$-hyponormal}
\author{Ra\'{u}l E. Curto}
\address{Department of Mathematics, The University of Iowa, Iowa City, Iowa
52242}
\email{rcurto@math.uiowa.edu}
\urladdr{http://www.math.uiowa.edu/\symbol{126}rcurto/}
\author{Sang Hoon Lee}
\address{Department of Mathematics, The University of Iowa, Iowa City, Iowa
52242}
\email{shlee@math.skku.ac.kr}
\urladdr{}
\thanks{The first named author was partially supported by NSF Grants
DMS-0099357 and DMS-0400741. \ The second named author was supported by the
Post-doctoral Fellowship Program of KOSEF }

\begin{abstract}
We give the first example of a quartically hyponormal unilateral weighted
shift which is not $3$-hyponormal.
\end{abstract}

\maketitle

\section{\label{Int}Introduction}

Let $\mathcal{H}$ and $\mathcal{K}$ be complex Hilbert spaces, let $\mathcal{%
L(H},\mathcal{K)}$ be the set of bounded linear operators from $\mathcal{H}$
to $\mathcal{K}$ and write $\mathcal{L(H)}:=\mathcal{L(H},\mathcal{H)}$. An
operator $T\in \mathcal{L(H)}$ is said to be \textit{normal} if $T^{\ast
}T=TT^{\ast }$, \textit{hyponormal} if $T^{\ast }T\geq TT^{\ast }$, and 
\textit{subnormal} if $T=N|_{\mathcal{H}}$, where $N$ is normal on some
Hilbert space $\mathcal{K}\supseteq \mathcal{H}$. \ If $T$ is subnormal then 
$T$ is also hyponormal. \ The Bram-Halmos criterion for subnormality states
that an operator $T$ is subnormal if and only if 
\begin{equation*}
\sum_{i,j}(T^{i}x_{j},T^{j}x_{i})\geq 0
\end{equation*}%
for all finite collections $x_{0},x_{1},\cdots ,x_{k}\in \mathcal{H}$ (\cite%
{Bra},\cite[II.1.9]{Con}). \ It is easy to see that this is equivalent to
the following positivity test: 
\begin{equation}
\left( 
\begin{array}{cccc}
I & T^{\ast } & \cdots & T^{\ast k} \\ 
T & T^{\ast }T & \cdots & T^{\ast k}T \\ 
\vdots & \vdots & \ddots & \vdots \\ 
T^{k} & T^{\ast }T^{k} & \cdots & T^{\ast k}T^{k}%
\end{array}%
\right) \geq 0\quad \text{(all $k\geq 1$)}.  \label{condition01}
\end{equation}%
Condition (\ref{condition01}) provides a measure of the gap between
hyponormality and subnormality. \ In fact, the positivity of (\ref%
{condition01}) for $k=1$ is equivalent to the hyponormality of $T$, while
subnormality requires the validity of (\ref{condition01}) for all $k$. \ 

Let $[A,B]:=AB-BA$ denote the commutator of two operators $A$ and $B$, and
define $T$ to be $k$-\textit{hyponormal} whenever the $k\times k$ operator
matrix 
\begin{equation}
M_{k}(T):=([T^{\ast j},T^{i}])_{i,j=1}^{k}  \label{condition02}
\end{equation}%
is positive. \ An application of the Choleski algorithm for operator
matrices shows that the positivity of (\ref{condition02}) is equivalent to
the positivity of the $(k+1)\times (k+1)$ operator matrix in (\ref%
{condition01}); the Bram-Halmos criterion can then be rephrased as saying
that $T$ is subnormal if and only if $T$ is $k$-hyponormal for every $k\geq
1 $ (\cite{CMX}).

Recall (\cite{Ath},\cite{CMX},\cite{CoS}) that $T\in \mathcal{L(H)}$ is said
to be \textit{weakly $k$-hyponormal} if 
\begin{equation*}
LS(T,T^{2},\cdots ,T^{k}):=\left\{ \sum_{j=1}^{k}\alpha _{j}T^{j}:\alpha
_{1},\cdots ,\alpha _{k}\in \mathbb{C}\right\}
\end{equation*}%
consists entirely of hyponormal operators, or equivalently, $M_{k}(T)$ is 
\textit{weakly positive}, cf. \cite{CMX}, i.e.,%
\begin{equation}
\left\langle M_{k}(T)\left( 
\begin{array}{c}
\lambda _{1}x \\ 
\vdots \\ 
\lambda _{k}x\ 
\end{array}%
\right) ,\left( 
\begin{array}{c}
\lambda _{1}x \\ 
\vdots \\ 
\lambda _{k}x\ 
\end{array}%
\right) \right\rangle \geq 0\quad \text{for all $x\in \mathcal{H}$ and $%
\lambda _{1},\cdots ,\lambda _{k}\in \mathbb{C}$}.  \label{weakhypo}
\end{equation}%
The operator $T$ is said to be \textit{quadratically hyponormal} when (\ref%
{weakhypo}) holds for $k=2$, and \textit{cubically hyponormal} (resp. 
\textit{quartically hyponormal}) when (\ref{weakhypo}) holds for $k=3$
(resp. $k=4$). \ Similarly, $T\in \mathcal{L(H)}$ is said to be \textit{%
polynomially hyponormal} if $p(T)$ is hyponormal for every polynomial $p\in 
\mathbb{C}[z]$. \ It is straightforward to verify that $k$-hyponormality
implies weak $k$-hyponormality, but the converse is not true in general. \
For unilateral weighted shifts, quadratic hyponormality is detected through
the analysis of an associated tridiagonal matrix, while cubic hyponormality
requires a pentadiagonal matrix\ (\cite{QHWS},\cite{JP}). \ The associated
nested determinants satisfy either a two-step recurring relation (in the
tridiagonal case) or a six-step recurring relation (in the pentadiagonal
case). \ 

The concrete calculation of the above mentioned nested determinants has
helped shed light on quadratic and cubic hyponormality. \ On the other hand,
quartic hyponormality requires heptadiagonal matrices, and a similar
multi-step recurring relation for the nested determinants is not known. \ As
a result, there is very little information available about quartic
hyponormality, and the notion has remained highly inscrutable.

In this paper, we present the first concrete example of a quartically
hyponormal operator which is not $3$-hyponormal. \ Although the result is
somewhat expected, and consistent with previous results in this area (e.g.,
for general operators polynomial hyponormality does not imply $2$%
-hyponormality \cite{CP1}, and for weighted shifts cubic hyponormality does
not imply $2$-hyponormality \cite{JP}), the techniques needed to prove it
are new. \ For instance, the proof of Theorem \ref{thm4}(i) includes a new
trick to compute a key determinant, while the proof of Theorem \ref{thm4}%
(ii) is based on a special rearrangement of the terms in the quadratic form $%
\Delta $ whose positivity ensures that the weighted shift is quartically
hyponormal.

Recall that given a bounded sequence of positive numbers $\alpha :\alpha
_{0},\alpha _{1},\cdots $ (called \textit{weights}), the \textit{%
(unilateral) weighted shift} $W_{\alpha }$ associated with $\alpha $ is the
operator on $\ell ^{2}(\mathbb{Z}_{+})$ defined by $W_{\alpha }e_{n}:=\alpha
_{n}e_{n+1}$ for all $n\geq 0$, where $\{e_{n}\}_{n=0}^{\infty }$ is the
canonical orthonormal basis for $\ell ^{2}$. It is straightforward to check
that $W_{\alpha }$ can never be normal, and that $W_{\alpha }$ is hyponormal
if and only if $\alpha _{n}\leq \alpha _{n+1}$ for all $n\geq 0$. \ The
moments of $\alpha $ are given as 
\begin{equation*}
\gamma _{k}\equiv \gamma _{k}(\alpha ):=\left\{ 
\begin{array}{cc}
1 & \text{if }k=0 \\ 
\alpha _{0}^{2}\cdots \alpha _{k-1}^{2} & \text{if }k>0.%
\end{array}%
\right.
\end{equation*}

We now recall a well known characterization of subnormality for
single-variable weighted shifts, due to C. Berger (cf. \cite[III.8.16]{Con}%
), and independently established by R. Gellar and L.J. Wallen \cite{GeWa}: $%
W_{\alpha }$ is subnormal if and only if there exists a probability measure $%
\xi $ (called the \textit{Berger measure} of $W_{\alpha }$) supported in $%
[0,\left\| W_{\alpha }\right\| ^{2}]$ such that $\gamma _{k}(\alpha
):=\alpha _{0}^{2}\cdots \alpha _{k-1}^{2}=\int t^{k}\;d\xi (t)\;\;(k\geq 1)$%
. \ If $W_{\alpha }$ is subnormal, and if for $h\geq 1$ we let $\mathcal{M}%
_{h}:=\bigvee \{e_{n}:n\geq h\}$ denote the invariant subspace obtained by
removing the first $h$ vectors in the canonical orthonormal basis of $\ell
^{2}(\mathbb{Z}_{+})$, then the Berger measure of $W_{\alpha }|_{\mathcal{M}%
_{h}}$ is $\frac{1}{\gamma _{h}}t^{h}d\xi (t)$.

The classes of (weakly) $k$-hyponormal operators have been studied in an
attempt to bridge the gap between subnormality and hyponormality (\cite%
{bridge}, \cite{QHWS}, \cite{OTAMP}, \cite{CuLe2}, \cite{CuLe3}, \cite{CLL2}%
, \cite{CMX}, \cite{DPY}, \cite{JP}, \cite{McCP}). \ The study of this gap
has been mostly successful at the level of $k$-hyponormality; for example,
for Toeplitz operators on the Hardy space of the unit circle, the gap is
described in (\cite{CLL2}). \ For weighted shifts, on the other hand,
positive results appear in \cite{QHWS} and \cite{CLL2}, although no concrete
example of a weighted shift which is polynomially hyponormal and not
subnormal has yet been found (the existence of such weighted shifts was
established in \cite{CP1} and \cite{CP2}). \ 

For weak $k$-hyponormality there nevertheless exist some partial results. \
For example, in \cite{QHWS}, the gap between $2$-hyponormality and quadratic
hyponormality was established, and in \cite{JP} weighted shifts which are
cubically hyponormal and not $2$-hyponormal were found. \ In this paper, we
give an example of a weighted shift which is weakly $4$-hyponormal but not $%
3 $-hyponormal.

\textit{Acknowledgment}. \ The authors are indebted to the referee for some
helpful suggestions. \ Many of the examples in this paper were obtained
using calculations with the software tool \textit{Mathematica \cite{Wol}.}

\section{Main Results}

We begin with an observation about quadratic hyponormality.

\begin{proposition}
$W_{\alpha }$ is quadratically hyponormal if and only if $W_{\alpha
}+sW_{\alpha }^{2}$ is hyponormal for all $s\geq 0.$
\end{proposition}

\begin{proof}
$(\Rightarrow )$ This implication is trivial.

$(\Leftarrow )$ Suppose $W_{\alpha }+sW_{\alpha }^{2}$ is hyponormal for all 
$s\geq 0$. $\ $We must show that $W_{\alpha }+cW_{\alpha }^{2}$ is
hyponormal for all $c\in \mathbb{C}.$ \ For $c\equiv se^{i\theta }\;(s>0)$,
there exists a unitary operator $U$ such that $UTU^{\ast }=e^{-i\theta }T.$
Then 
\begin{eqnarray*}
U(T+cT^{2})U^{\ast } &=&UTU^{\ast }+cUT^{2}U^{\ast } \\
&=&UTU^{\ast }+c(UTU^{\ast })^{2} \\
&=&e^{-i\theta }T+se^{i\theta }\cdot e^{-2i\theta }T^{2} \\
&=&e^{-i\theta }(T+sT^{2})
\end{eqnarray*}%
is hyponormal. Therefore, $T+cT^{2}$ is hyponormal.
\end{proof}

\begin{lemma}
The following statements are equivalent.\newline
(i) $W_{\alpha }$ is quartically hyponormal; \newline
(ii) For each $x\equiv \{x_{n}\}_{n=0}^{\infty }\in \ell ^{2}$, we have $%
(\left\langle [W_{\alpha }^{\ast j},W_{\alpha }^{i}]x,x\right\rangle
)_{i,j=1}^{4}\geq 0;$ \newline
(iii) For each $a,b,c\in \mathbb{C}$ and $x\equiv \{x_{n}\}_{n=0}^{\infty
}\in \ell ^{2}$, 
\begin{multline*}
\Delta :=|c|^{2}r_{0}|x_{0}|^{2}+\left\langle \Theta _{1}\left( 
\begin{array}{c}
\overline{b}x_{0} \\ 
\overline{c}x_{1}%
\end{array}%
\right) ,\left( 
\begin{array}{c}
\overline{b}x_{0} \\ 
\overline{c}x_{1}%
\end{array}%
\right) \right\rangle +\left\langle \Theta _{2}\left( 
\begin{array}{c}
\overline{a}x_{0} \\ 
\overline{b}x_{1} \\ 
\overline{c}x_{2}%
\end{array}%
\right) ,\left( 
\begin{array}{c}
\overline{a}x_{0} \\ 
\overline{b}x_{1} \\ 
\overline{c}x_{2}%
\end{array}%
\right) \right\rangle \\
+\sum_{i=0}^{\infty }\left\langle \Delta _{i}\left( 
\begin{array}{c}
x_{i} \\ 
\overline{a}x_{i+1} \\ 
\overline{b}x_{i+2} \\ 
\overline{c}x_{i+3}%
\end{array}%
\right) ,\left( 
\begin{array}{c}
x_{i} \\ 
\overline{a}x_{i+1} \\ 
\overline{b}x_{i+2} \\ 
\overline{c}x_{i+3}%
\end{array}%
\right) \right\rangle \geq 0.
\end{multline*}%
Here 
\begin{equation*}
\Theta _{1}:=\left( 
\begin{array}{cc}
p_{0} & \sqrt{g_{0}} \\ 
\sqrt{g_{0}} & r_{1}%
\end{array}%
\right) ,\;\;\Theta _{2}:=\left( 
\begin{array}{ccc}
v_{0} & \sqrt{t_{0}} & \sqrt{f_{0}} \\ 
\sqrt{t_{0}} & p_{1} & \sqrt{g_{1}} \\ 
\sqrt{f_{0}} & \sqrt{g_{1}} & r_{2}%
\end{array}%
\right) ,
\end{equation*}%
\begin{equation*}
\Delta _{i}:=\left( 
\begin{array}{cccc}
u_{i} & \sqrt{w_{i}} & \sqrt{s_{i}} & \sqrt{q_{i}} \\ 
\sqrt{w_{i}} & v_{i+1} & \sqrt{t_{i+1}} & \sqrt{f_{i+1}} \\ 
\sqrt{s_{i}} & \sqrt{t_{i+1}} & p_{i+2} & \sqrt{g_{i+2}} \\ 
\sqrt{q_{i}} & \sqrt{f_{i+1}} & \sqrt{g_{i+2}} & r_{i+3}%
\end{array}%
\right) \;(i\geq 0),
\end{equation*}%
where 
\begin{eqnarray*}
u_{i} &:&=\alpha _{i}^{2}-\alpha _{i-1}^{2} \\
w_{i} &:&=\alpha _{i}^{2}(\alpha _{i+1}^{2}-\alpha _{i-1}^{2})^{2} \\
v_{i} &:&=\alpha _{i}^{2}\alpha _{i+1}^{2}-\alpha _{i-1}^{2}\alpha _{i-2}^{2}
\\
s_{i} &:&=\alpha _{i}^{2}\alpha _{i+1}^{2}(\alpha _{i+2}^{2}-\alpha
_{i-1}^{2})^{2} \\
t_{i} &:&=\alpha _{i}^{2}(\alpha _{i+1}^{2}\alpha _{i+2}^{2}-\alpha
_{i-1}^{2}\alpha _{i-2}^{2})^{2} \\
p_{i} &:&=\alpha _{i}^{2}\alpha _{i+1}^{2}\alpha _{i+2}^{2}-\alpha
_{i-1}^{2}\alpha _{i-2}^{2}\alpha _{i-3}^{2} \\
q_{i} &:&=\alpha _{i}^{2}\alpha _{i+1}^{2}\alpha _{i+2}^{2}(\alpha
_{i+3}^{2}-\alpha _{i-1}^{2})^{2} \\
f_{i} &:&=\alpha _{i}^{2}\alpha _{i+1}^{2}(\alpha _{i+2}^{2}\alpha
_{i+3}^{2}-\alpha _{i-1}^{2}\alpha _{i-2}^{2})^{2} \\
g_{i} &:&=\alpha _{i}^{2}(\alpha _{i+1}^{2}\alpha _{i+2}^{2}\alpha
_{i+3}^{2}-\alpha _{i-1}^{2}\alpha _{i-2}^{2}\alpha _{i-3}^{2})^{2} \\
r_{i} &:&=\alpha _{i}^{2}\alpha _{i+1}^{2}\alpha _{i+2}^{2}\alpha
_{i+3}^{2}-\alpha _{i-1}^{2}\alpha _{i-2}^{2}\alpha _{i-3}^{2}\alpha
_{i-4}^{2}.
\end{eqnarray*}%
(As usual, we let $\alpha _{-1}=\alpha _{-2}=\alpha _{-3}=\alpha _{-4}=0$.)
\end{lemma}

\begin{proof}
This is a straightforward computation.
\end{proof}

\begin{remark}
Observe that $W_{\alpha }$ is $4$-hyponormal if and only if $\Theta _{2}\
\geq 0$ and $\Delta _{i}\geq 0$ for all $i\geq 0.$
\end{remark}

We now give an example of a weighted shifts which is quartically hyponormal
but not $3$-hyponormal.

\begin{theorem}
\label{thm4}For $x>0$, let $W_{\alpha (x)}$ be the unilateral weighted shift
with weight sequence given by $\alpha _{0}:=\sqrt{x},\quad \alpha _{n}:=%
\sqrt{\frac{n+2}{n+3}}\quad (n\geq 1)$. Then\newline
(i) $W_{\alpha (x)}$ is $k$-hyponormal if and only if $0<x\leq \frac{%
2(k+1)^{2}(k+2)^{2}}{3k(k+3)(k^{2}+3k+4)}=:H_{k}$;\newline
(ii) If $0<x\leq \frac{667}{990}$, then $W_{\alpha (x)}$ is quartically
hyponormal.
\end{theorem}

\begin{proof}
(i) By \cite[Theorem 4(d)]{QHWS}, we know that $W_{\alpha (x)}$ is $k$%
-hyponormal if and only if 
\begin{equation*}
A(n;k):=\left( 
\begin{array}{cccc}
\gamma _{n} & \gamma _{n+1} & \cdots & \gamma _{n+k} \\ 
\gamma _{n+1} & \gamma _{n+2} & \cdots & \gamma _{n+k+1} \\ 
\vdots & \vdots & \ddots & \vdots \\ 
\gamma _{n+k} & \gamma _{n+k+1} & \cdots & \gamma _{n+2k}%
\end{array}%
\right) \geq 0\quad \text{(all}\quad n\geq 0\text{).}
\end{equation*}%
Since $W_{\alpha (x)}$ has a Bergman tail, it is enough to check at $n=0$. \
In this case, $A(0;k)\geq 0$ is equivalent to 
\begin{equation*}
\det \left( 
\begin{array}{ccccc}
\frac{1}{3x} & \frac{1}{3} & \frac{1}{4} & \cdots & \frac{1}{k+2} \\ 
\frac{1}{3} & \frac{1}{4} & \frac{1}{5} & \cdots & \frac{1}{k+3} \\ 
\vdots & \vdots & \vdots & \ddots & \vdots \\ 
\frac{1}{k+2} & \frac{1}{k+3} & \frac{1}{k+4} & \cdots & \frac{1}{2k+2}%
\end{array}%
\right) \geq 0.
\end{equation*}%
Let 
\begin{equation*}
A:=\left( 
\begin{array}{ccccc}
\frac{1}{3x} & \frac{1}{3} & \frac{1}{4} & \cdots & \frac{1}{k+2} \\ 
\frac{1}{3} & \frac{1}{4} & \frac{1}{5} & \cdots & \frac{1}{k+3} \\ 
\vdots & \vdots & \vdots & \ddots & \vdots \\ 
\frac{1}{k+2} & \frac{1}{k+3} & \frac{1}{k+4} & \cdots & \frac{1}{2k+2}%
\end{array}%
\right) ,\;B:=\left( 
\begin{array}{ccccc}
\frac{1}{2} & \frac{1}{3} & \frac{1}{4} & \cdots & \frac{1}{k+2} \\ 
\frac{1}{3} & \frac{1}{4} & \frac{1}{5} & \cdots & \frac{1}{k+3} \\ 
\vdots & \vdots & \vdots & \ddots & \vdots \\ 
\frac{1}{k+2} & \frac{1}{k+3} & \frac{1}{k+4} & \cdots & \frac{1}{2k+2}%
\end{array}%
\right)
\end{equation*}%
and%
\begin{equation*}
C:=\left( 
\begin{array}{cccc}
\frac{1}{4} & \frac{1}{5} & \cdots & \frac{1}{k+3} \\ 
\frac{1}{5} & \frac{1}{6} & \cdots & \frac{1}{k+4} \\ 
\vdots & \vdots & \ddots & \vdots \\ 
\frac{1}{k+3} & \frac{1}{k+4} & \cdots & \frac{1}{2k+2}%
\end{array}%
\right) .
\end{equation*}%
Expanding the determinants of $A$ and $B$ by the first row, we have 
\begin{eqnarray*}
\det A &=&\frac{1}{3x}\det C+Q \\
\det B &=&\frac{1}{2}\det C+Q,
\end{eqnarray*}%
so that 
\begin{eqnarray*}
\det A &=&\frac{1}{3x}\det C+\det B-\frac{1}{2}\det C \\
&=&\frac{2-3x}{6x}\det C+\det B.
\end{eqnarray*}

Now, for $H:=(h_{ij})_{i,j=1}^{n},h_{ij}:=(p+i+j-1)^{-1}$ and $p\geq 0$,
recall that 
\begin{equation*}
\det H=(1!2!\cdots (n-1)!)^{2}\frac{\Gamma (p+1)\Gamma (p+2)\cdots \Gamma
(p+n)}{\Gamma (n+p+1)\Gamma (n+p+2)\cdots \Gamma (2n+p)}.
\end{equation*}%
Thus, 
\begin{equation*}
\det A=\frac{(1!2!\cdots (k-1)!)^{2}\Gamma (4)\Gamma (5)\cdots \Gamma (k+2)}{%
\Gamma (k+4)\Gamma (k+5)\cdots \Gamma (2k+3)}[(\frac{2-3x}{6x})\Gamma (k+3)+%
\frac{(k!)^{2}\Gamma (2)\Gamma (3)}{\Gamma (k+3)}].
\end{equation*}%
Therefore, $\det A\geq 0$ if and only if $0<x\leq \frac{2(k+1)^{2}(k+2)^{2}}{%
3k(k+3)(k^{2}+3k+4)}$, as desired.

(ii) By a direct computation, we have 
\begin{equation*}
\Theta _{1}=\left( 
\begin{array}{cc}
\frac{3}{5}x & \frac{\sqrt{x}}{2} \\ 
\frac{\sqrt{x}}{2} & \frac{3}{7}%
\end{array}%
\right) ,\quad \Theta _{2}=\left( 
\begin{array}{ccc}
\frac{3}{4}x & \frac{3}{5}\sqrt{x} & \sqrt{\frac{x}{3}} \\ 
\frac{3}{5}\sqrt{x} & \frac{1}{2} & \frac{2\,\sqrt{3}}{7} \\ 
\sqrt{\frac{x}{3}} & \frac{2\,\sqrt{3}}{7} & \frac{1}{2}%
\end{array}%
\right) ,
\end{equation*}%
\begin{equation*}
\Delta _{0}=\left( 
\begin{array}{cccc}
x & \frac{3}{4}\sqrt{x} & \frac{2}{5}\sqrt{3x} & \frac{1}{2}\sqrt{\frac{5}{3}%
x} \\ 
\frac{3}{4}\sqrt{x} & \frac{3}{5} & \frac{1}{\sqrt{3}} & \frac{\sqrt{15}}{7}
\\ 
\frac{2}{5}\sqrt{3x} & \frac{1}{\sqrt{3}} & \frac{4}{7} & \frac{\sqrt{5}}{4}
\\ 
\frac{1}{2}\sqrt{\frac{5}{3}x} & \frac{\sqrt{15}}{7} & \frac{\sqrt{5}}{4} & 
\frac{5}{9}%
\end{array}%
\right) ,
\end{equation*}%
\begin{equation*}
\Delta _{1}=\left( 
\begin{array}{cccc}
\frac{3}{4}-x & \frac{\sqrt{3}}{2}(\frac{4}{5}-x) & \sqrt{\frac{3}{5}}(\frac{%
5}{6}-x) & \frac{1}{\sqrt{2}}(\frac{6}{7}-x) \\ 
\frac{\sqrt{3}}{2}(\frac{4}{5}-x) & \frac{2}{3}-\frac{3}{4}x & \frac{2}{%
\sqrt{5}}(\frac{5}{7}-\frac{3}{4}x) & \frac{\sqrt{3}}{2\,\sqrt{2}}(1-x) \\ 
\sqrt{\frac{3}{5}}(\frac{5}{6}-x) & \frac{2}{\sqrt{5}}(\frac{5}{7}-\frac{3}{4%
}x) & \frac{5}{8}-\frac{3}{5}x & \sqrt{\frac{5}{6}}(\frac{2}{3}-\frac{3}{5}x)
\\ 
\frac{1}{\sqrt{2}}(\frac{6}{7}-x) & \frac{\sqrt{3}}{2\,\sqrt{2}}(1-x) & 
\sqrt{\frac{5}{6}}(\frac{2}{3}-\frac{3}{5}x) & \frac{3}{5}-\frac{x}{2}%
\end{array}%
\right) ,
\end{equation*}%
\begin{equation*}
\Delta _{2}=\left( 
\begin{array}{cccc}
\frac{1}{20} & \frac{1}{6\sqrt{5}} & \frac{\sqrt{\frac{3}{2}}}{14} & \frac{1%
}{4\sqrt{7}} \\ 
\frac{1}{6\sqrt{5}} & \frac{4}{35} & \frac{\sqrt{\frac{3}{10}}}{4} & \frac{8%
}{9\sqrt{35}} \\ 
\frac{\sqrt{\frac{3}{2}}}{14} & \frac{\sqrt{\frac{3}{10}}}{4} & \frac{1}{6}
& \frac{\sqrt{\frac{6}{7}}}{5} \\ 
\frac{1}{4\sqrt{7}} & \frac{8}{9\sqrt{35}} & \frac{\sqrt{\frac{6}{7}}}{5} & 
\frac{16}{77}%
\end{array}%
\right) ,\;\;\Delta _{3}=\left( 
\begin{array}{cccc}
\frac{1}{30} & \frac{\sqrt{\frac{2}{15}}}{7} & \frac{3}{8\sqrt{35}} & \frac{%
\sqrt{\frac{2}{5}}}{9} \\ 
\frac{\sqrt{\frac{2}{15}}}{7} & \frac{1}{12} & \frac{\sqrt{\frac{2}{21}}}{3}
& \frac{1}{5\sqrt{3}} \\ 
\frac{3}{8\sqrt{35}} & \frac{\sqrt{\frac{2}{21}}}{3} & \frac{9}{70} & \frac{3%
\sqrt{\frac{2}{7}}}{11} \\ 
\frac{\sqrt{\frac{2}{5}}}{9} & \frac{1}{5\sqrt{3}} & \frac{3\sqrt{\frac{2}{7}%
}}{11} & \frac{1}{6}%
\end{array}%
\right) ,
\end{equation*}%
and 
\begin{equation*}
\Delta _{4}=\left( 
\begin{array}{cccc}
\frac{1}{42} & \frac{1}{4\sqrt{42}} & \frac{1}{12\sqrt{3}} & \frac{\sqrt{%
\frac{2}{3}}}{15} \\ 
\frac{1}{4\sqrt{42}} & \frac{4}{63} & \frac{3}{10\sqrt{14}} & \frac{8}{33%
\sqrt{7}} \\ 
\frac{1}{12\sqrt{3}} & \frac{3}{10\sqrt{14}} & \frac{9}{88} & \frac{1}{6%
\sqrt{2}} \\ 
\frac{\sqrt{\frac{2}{3}}}{15} & \frac{8}{33\sqrt{7}} & \frac{1}{6\sqrt{2}} & 
\frac{16}{117}%
\end{array}%
\right) .
\end{equation*}%
Note that $\Delta _{n}\geq 0$ for all $n\geq 5$, and that all above matrices
are positive except possibly for $\Delta _{1}$. \ For the positivity of $%
\Delta $, we minimize the positivity of $\Theta _{1},\Theta _{2},\Delta
_{0},\Delta _{2},\Delta _{3},\Delta _{4}$, that is, we replace $\Theta
_{1},\Theta _{2},\Delta _{0},\Delta _{2},\Delta _{3},\Delta _{4}$ by $\Theta
_{1}^{\prime },\Theta _{2}^{\prime },\Delta _{0}^{\prime },\Delta
_{2}^{\prime },$ $\Delta _{3}^{\prime },\Delta _{4}^{\prime }$, with $%
\operatorname{rank}\Theta _{j}^{\prime }=j$ for $j=1,2$, and $\operatorname{rank}%
\Delta _{j}^{\prime }=3$ for $j=0,2,3,4$, respectively, where 
\begin{equation*}
\Theta _{1}^{\prime }:=\left( 
\begin{array}{cc}
\cdot & \cdot \\ 
\cdot & \frac{5}{12}%
\end{array}%
\right) ,\quad \Theta _{2}^{\prime }:=\left( 
\begin{array}{ccc}
\cdot & \cdot & \cdot \\ 
\cdot & \frac{612}{1225} & \cdot \\ 
\cdot & \cdot & \cdot%
\end{array}%
\right) ,
\end{equation*}%
$\newline
$%
\begin{equation*}
\Delta _{0}^{\prime }:=\left( 
\begin{array}{cccc}
\cdot & \cdot & \cdot & \cdot \\ 
\cdot & \frac{1411}{2352} & \cdot & \cdot \\ 
\cdot & \cdot & \cdot & \cdot \\ 
\cdot & \cdot & \cdot & \cdot%
\end{array}%
\right) ,\quad \Delta _{2}^{\prime }:=\left( 
\begin{array}{cccc}
\frac{627}{12544} & \cdot & \cdot & \cdot \\ 
\cdot & \cdot & \cdot & \cdot \\ 
\cdot & \cdot & \cdot & \cdot \\ 
\cdot & \cdot & \cdot & \cdot%
\end{array}%
\right) ,
\end{equation*}%
\begin{equation*}
\Delta _{3}^{\prime }:=\left( 
\begin{array}{cccc}
\frac{1411}{42336} & \cdot & \cdot & \cdot \\ 
\cdot & \cdot & \cdot & \cdot \\ 
\cdot & \cdot & \cdot & \cdot \\ 
\cdot & \cdot & \cdot & \cdot%
\end{array}%
\right) \;\text{and }\Delta _{4}^{\prime }:=\left( 
\begin{array}{cccc}
\frac{2057}{86400} & \cdot & \cdot & \cdot \\ 
\cdot & \cdot & \cdot & \cdot \\ 
\cdot & \cdot & \cdot & \cdot \\ 
\cdot & \cdot & \cdot & \cdot%
\end{array}%
\right) .
\end{equation*}%
(In all of the above expressions, the symbol ``$\cdot $'' denotes an entry
that remains unchanged.) \ Let $\widetilde{\Delta }(a,b,c):=$ 
\begin{equation*}
\left( 
\begin{array}{cccc}
A & \frac{\sqrt{3}}{2}(\frac{4}{5}-x)\bar{a} & \sqrt{\frac{3}{5}}(\frac{5}{6}%
-x)\bar{b} & \frac{1}{\sqrt{2}}(\frac{6}{7}-x)\bar{c} \\ 
\frac{\sqrt{3}}{2}(\frac{4}{5}-x)a & B & \frac{2}{\sqrt{5}}(\frac{5}{7}-%
\frac{3}{4}x)a\bar{b} & \frac{\sqrt{3}}{2\,\sqrt{2}}(1-x)a\bar{c} \\ 
\sqrt{\frac{3}{5}}(\frac{5}{6}-x)b & \frac{2}{\sqrt{5}}(\frac{5}{7}-\frac{3}{%
4}x)\bar{a}b & (\frac{5}{8}-\frac{3}{5}x)|b|^{2}+\frac{1}{211680} & \sqrt{%
\frac{5}{6}}(\frac{2}{3}-\frac{3}{5}x)b\bar{c} \\ 
\frac{1}{\sqrt{2}}(\frac{6}{7}-x)c & \frac{\sqrt{3}}{2\,\sqrt{2}}(1-x)\bar{a}%
c & \sqrt{\frac{5}{6}}(\frac{2}{3}-\frac{3}{5}x)\bar{b}c & (\frac{3}{5}-%
\frac{x}{2})|c|^{2}+\frac{1}{604800}%
\end{array}%
\right) ,
\end{equation*}%
where $A:=\frac{3}{4}-x+\frac{|a|^{2}}{11760}+\frac{|b|^{2}}{2450}+\frac{%
|c|^{2}}{84}$ and $B:=(\frac{2}{3}-\frac{3}{4}x)|a|^{2}+\frac{1}{62720}$. \
If $\widetilde{\Delta }(a,b,c)\geq 0$ for all $a,b,c\in \mathbb{C}$, then $%
\Delta \geq 0$ and hence, by Lemma 2, $W_{\alpha (x)}$ is quartically
hyponormal. \ Note that in each nested determinant of $\widetilde{\Delta }%
(a,b,c)$ the parameters $a,b,c$ occur in modulus square form. \ So, using
the Nested Determinants Test, we can easily see that if $x\leq \frac{22580899%
}{33531912}=:\xi $ then every coefficient in $\det \widetilde{\Delta }%
(a,b,c) $ is positive. \ Therefore, $\widetilde{\Delta }(a,b,c)\geq 0$ for
all $a,b,c\in \mathbb{C}.$ \ Note that $H_{2}=\frac{24}{35}$ and $H_{3}=%
\frac{200}{297}$. \ Observe that $H_{3}<\xi <H_{2}.$ \ Moreover, we can show
that $\widetilde{\Delta }(a,b,c)\geq 0$ for $x=\frac{667}{990}>\xi $, again
using the Nested Determinants Test.
\end{proof}

\begin{corollary}
(i) If $\frac{200}{297}<x\leq \frac{667}{990}$, then $W_{\alpha (x)}$ is
quartically hyponormal but not $3$-hyponormal. \newline
(ii) If $W_{\alpha (x)}$ is $3$-hyponormal then it is also quartically
hyponormal.
\end{corollary}

\end{document}